# SOME UNIFIED INTEGRALS ASSOCIATED WITH GENERALIZED STRUVE FUNCTION

**D.L. Suthar[1], S.D. Purohit[2] and K.S. Nisar[3]**

This paper is devoted for the study of a new generalization of Struve function type. In this paper , We establish four new integral formulas involving the Galué type Struve function, which are express in term of the generalized (Wright) hypergeometric functions. The result established here are general in nature and are likely to find useful in applied problem of Sciences, engineering and technology.

**Keywords**: Galué type Struve function, Gamma Function, Generalized Wright function, Oberhettinger integral formula, Lavoie-Trottier integral formula.

**AMS Subject Classification:** 33B20, 33C20, 33B15, 33C05.

## 1. INTRODUCTION

In recent year, the interest reader in the fractional calculus operators involving various special functions have found significant importance and applications in modeling of relevant systems in various fields science, engineering and technology,

Recently, Nisar et al. [12], defined as following generalized form of Struve function named as generalized Galué type Struve function (GTSF) as:

$$_a w_{p,b,c,\xi}^{\lambda,\mu}(z) = \sum_{k=0}^{\infty} \frac{(-c)^k}{\Gamma(\lambda k + \mu)\Gamma\left(ak + \frac{p}{\xi} + \frac{b+2}{2}\right)}\left(\frac{z}{2}\right)^{2k+p+1}, \quad a \in \mathbb{N} \; p,b,c \in \mathbb{C} \tag{1.1}$$

where $\lambda > 0, \; \xi > 0$ and $\mu$ is an arbitrary parameter.

For the definition of the struve function and its more generalization, the interested reader may refer to the many papers ( Bhow-mick [3, 4], Kanth [6], Singh [17, 18], Nisar and Atangana [10], Singh [19].

**Special cases**

When $\lambda = a = 1, \; \mu = 3/2$ and $\xi = 1$ in equation (1.1), it reduces to generalization of Struve function which is defined by Orhan and Yagmur [15, 16 ].

$$H_{p,b,c}(z) = \sum_{k=0}^{\infty} \frac{(-c)^k}{\Gamma\left(k + \frac{3}{2}\right)\Gamma\left(k + p + \frac{b+2}{2}\right)}\left(\frac{z}{2}\right)^{2k+p+1}, \quad p,b,c \in \mathbb{C} \tag{1.2}$$

Details related to the function $H_{p,b,c}(z)$ and its particular cases can be seen in Barics [1, 2], Mondal and Swaminathan [9], Mondal and Nisar [8]. Nisar et al. [11, 13].

For our purpose, we recall the generaliged Wright hypergeometric function $_p\psi_q(z)$ (see, for detail, Srivastava and Karson [21]), for $z \in \mathbb{C}$ complex, $a_i, \; b_j \in \mathbb{C}$ and $\alpha_i, \; \beta_j \in \Re$ , where $(\alpha_i, \; \beta_j \neq 0;$ $i = 1,2,...,p; j = 1, 2,...,q)$ , is defined as below:

$$_p\psi_q(z) = {}_p\psi_q\left[\begin{matrix}(a_i, \alpha_i)_{1,\; p} \\ (b_j, \beta_j)_{1,\; q}\end{matrix}\bigg| z\right] = \sum_{k=0}^{\infty} \frac{\prod_{i=1}^{p} \Gamma(a_i + \alpha_i k)\, z^k}{\prod_{j=1}^{q} \Gamma(b_j + \beta_j k)\, k!} \tag{1.3}$$

Introduced by Wright [22], the generalized Wright function (1.12) and proved several theorems on the asymptotic expansion of $_p\psi_q(z)$ (for instance, see [22, 23, 24]) for all values of the argument $z$, under the condition:

$$\sum_{j=1}^{q}\beta_j - \sum_{i=1}^{p}\alpha_i > -1. \tag{1.4}$$

It is noted that the generalized (Wright) hypergeometric function $_p\psi_q$ in (1.3) whose asymptotic expansion was investigated by Fox [5] and Wright is an interesting further generalization of the generalized hypergeometric series (1.6) as follow:

$$_p\Psi_q\left[\begin{matrix}(a_1,1),...,(a_p,1)\\(b_1,1),...,(b_q,1)\end{matrix}\; z\right] = \frac{\prod_{j=1}^{p}\Gamma(\alpha_j)}{\prod_{j=1}^{q}\Gamma(\beta_j)}\; _pF_q\left[\begin{matrix}\alpha_1,...,\alpha_p;\\\beta_1,...,\beta_q;\end{matrix}\; z\right], \tag{1.5}$$

where $_pF_q$ is the generalized hypergeometric series defined by (see [20, Section 1.5])

$$_pF_q\left[\begin{matrix}\alpha_1,...,\alpha_p;\\\beta_1,...,\beta_q;\end{matrix}\; z\right] = \sum_{n=0}^{\infty}\frac{(\alpha_1)_n\cdots(\alpha_p)_n}{(\beta_1)_n\cdots(\beta_q)_n}\frac{z^n}{n!} \quad = {_pF_q(\alpha_1,...,\alpha_p;\beta_1,...,\beta_q;z)}, \tag{1.6}$$

where $(\lambda)_n$ is the Pochhammer symbol defined (for $\lambda \in \mathbb{C}$) by (see [20, p.2 and p.4-6]):

$$(\lambda)_n = \begin{cases}1 & (n=0)\\\lambda(\lambda+1)\dots(\lambda+n-1) & (n \in N = \{1,2,3,...\})\end{cases} = \frac{\Gamma(\lambda+n)}{\Gamma(\lambda)} \quad (\lambda \in \mathbb{C}\setminus\mathbb{Z}_0) \tag{1.7}$$

and $\mathbb{Z}_0$ denotes the set of nonpositive integers.

For our present investigation, we also need to recall the following Oberhettinger's integral formula [14]:

$$\int_0^\infty x^{\mu-1}\left(x+a+\sqrt{x^2+2ax}\right)^{-\lambda}dx = 2\lambda a^{-\lambda}\left(\frac{a}{2}\right)^\mu\frac{\Gamma(2\mu)\Gamma(\lambda-\mu)}{\Gamma(1+\lambda+\mu)}, \tag{1.8}$$

provided $0 < \Re(\mu) < \Re(\lambda)$.

Also we recall Lavoie-Trottier integral formula [7] for our present study

$$\int_0^1 x^{\alpha-1}(1-x)^{2\beta-1}\left(1-\frac{x}{3}\right)^{2\alpha-1}\left(1-\frac{x}{4}\right)^{\beta-1}dx = \left(\frac{2}{3}\right)^{2\alpha}\frac{\Gamma(\alpha)\Gamma(\beta)}{\Gamma(\alpha+\beta)}, \tag{1.9}$$

provided $\Re(\alpha)>0, \;\Re(\beta)>0.$

## 2. Main Results

We establish four generalized integral formulas, which are expressed in terms of the generalized (Wright) hypergeometric functions (1.3), by inserting the generalized Galué type Struve function (1.1) with suitable arguments into the integrand of the integral (1.8).

**Theorem 1.** Let $x>0$ ; $a\in \mathbb{N}$ $\lambda, p,b,c\in \mathbb{C}$; $\upsilon>0$ and $\delta$ is an arbitrary parameter be such that $0 < \Re(\mu) < \Re(\lambda+p+1)$ then there hold the following results:

$$\int_0^\infty x^{\mu-1}\left(x+a+\sqrt{x^2+2ax}\right)^{-\lambda}{_aw_{p,b,c,\xi}^{\upsilon,\delta}}\left(\frac{y}{x+a+\sqrt{x^2+2ax}}\right)dx$$

$$= 2^{-\mu-p}\, a^{\mu-\lambda-p-1}\, y^{p+1}\, \Gamma(2\mu)\; {}_3\psi_4 \left[ \begin{array}{c} (\lambda+p+2,\,2),\;\; (\lambda-\mu+p+1,\,2),\,(1,1); \\[2mm] (\delta,\upsilon),\;\left(\dfrac{p}{\xi}+\dfrac{b}{2}+1,\;a\right),\;(\lambda+p+1,\;2),\,(\lambda+\mu+p+2,\;2) \end{array} \;\middle|\; \dfrac{-cy^2}{4a^2} \right] \qquad (2.1)$$

***Theorem 2.*** *Let* $x>0$ *;* $a\in\mathbb{N}$ *,* $\lambda,p,b,c\in\mathbb{C}$ *;* $\upsilon>0$ *and* $\delta$ *is an arbitrary parameter be such that* $0<\Re(\mu)<\Re(\lambda+p+1)$ *then there hold the following results:*

$$\int_0^\infty x^{\mu-1}\left(x+a+\sqrt{x^2+2ax}\right)^{-\lambda}{}_a w_{p,b,c,\xi}^{\upsilon,\delta}\left(\dfrac{xy}{x+a+\sqrt{x^2+2ax}}\right)dx$$

$$= 2^{-\mu-2p}\, a^{\mu-\lambda-1}\, y^{p+1}\, \Gamma(\lambda-\mu+1)\; {}_3\psi_4 \left[ \begin{array}{c} (\lambda+p+2,\,2),\;\; (2\mu+2p,\,4),\,(1,1); \\[2mm] (\delta,\upsilon),\;\left(\dfrac{p}{\xi}+\dfrac{b}{2}+1,\;a\right),\;(\lambda+p+1,\;2),\,(\lambda+\mu+2p+2,\,4); \end{array} \;\middle|\; \dfrac{-cy^2}{4} \right] \!(2.2)$$

***Proof.*** By making use of (1.1) in the ingrand of (2.1) and then interchanging the order of integral sign and summation, which is verified by uniform convergence of the involved series under the given conditions, we get

$$\int_0^\infty x^{\mu-1}\left(x+a+\sqrt{x^2+2ax}\right)^{-\lambda}{}_a w_{p,b,c,\xi}^{\upsilon,\delta}\left(\dfrac{y}{x+a+\sqrt{x^2+2ax}}\right)dx$$

$$= \sum_{k=0}^\infty \dfrac{(-c)^k}{\Gamma(\upsilon k+\delta)\Gamma\left(ak+\dfrac{p}{\xi}+\dfrac{b+2}{2}\right)}\left(\dfrac{y}{2}\right)^{2k+p+1}\int_0^\infty x^{\mu-1}\left(x+a+\sqrt{x^2+2ax}\right)^{-(\lambda+2k+p+1)}dx \qquad (2.3)$$

we can apply the integral formula (1.8) to the integral in (2.3) and obtain the following expression:

$$= \sum_{k=0}^\infty \dfrac{(-c)^k}{\Gamma(\upsilon k+\delta)\Gamma\left(ak+\dfrac{p}{\xi}+\dfrac{b+2}{2}\right)}\left(\dfrac{y}{2}\right)^{2k+p+1}\left(\dfrac{a}{2}\right)^\mu 2(\lambda+p+1+2k)\,a^{-(\lambda+p+1+2k)}\dfrac{\Gamma(2\mu)\Gamma(\lambda-\mu+p+1+2k)}{\Gamma(\lambda+\mu+p+2+2k)}$$

$$= 2^{-\mu-p}\, a^{\mu-\lambda-p-1}\, y^{p+1}\, \Gamma(2\mu)\sum_{k=0}^\infty \dfrac{(-c)^k\,\Gamma(\lambda+p+2+2k)\,\Gamma(\lambda-\mu+p+1+2k)}{\Gamma(\upsilon k+\delta)\Gamma\left(ak+\dfrac{p}{\xi}+\dfrac{b+2}{2}\right)\Gamma(\lambda+p+1+2k)\Gamma(\lambda+\mu+p+2+2k)}\left(\dfrac{y}{2a}\right)^{2k}$$

In accordance with the definition of (1.3), we obatain the result (2.1). This completes the proof of the theorem.

By similar manner as in proof of Theorem 1, we can prove the integral formula (2.2).

***Theorem 3.*** *Let* $x>0$ *;* $a\in\mathbb{N}$ *,* $\lambda,p,b,c\in\mathbb{C}$ *;* $\upsilon>0$ *and* $\delta$ *is an arbitrary parameter be such that* $\Re(\alpha)>0$ *,* $\Re(\beta+p+1+2k)>0$ *then there hold the following results:*

$$\int_0^1 x^{\alpha-1}(1-x)^{2\beta-1}\left(1-\dfrac{x}{3}\right)^{2\alpha-1}\left(1-\dfrac{x}{4}\right)^{\beta-1}{}_a w_{p,b,c,\xi}^{\upsilon,\delta}\left(y\left(1-\dfrac{x}{4}\right)(1-x)^2\right)dx$$

$$= \left(\frac{2}{3}\right)^{2\alpha} \left(\frac{y}{2}\right)^{p+1} \Gamma(2\alpha) \,_2\psi_3 \left[ \begin{array}{c} (\beta+p+1,\,2),\,(1,1); \\ \\ (\delta,\,\upsilon),\,\left(\frac{p}{\xi}+\frac{b}{2}+1,\,a\right),\,(2\alpha+\beta+p+1,2); \end{array} \quad \frac{-cy^2}{4} \right] \tag{2.4}$$

**Theorem 4.** *Let* $x>0$ ; $a \in \mathbb{N}$ $\lambda, p,b,c \in \mathbb{C}$; $\upsilon>0$ *and* $\delta$ *is an arbitrary parameter be such that* $\Re(\beta)>0$, $\Re(\alpha+p+1+2k)>0$ *then there hold the following results:*

$$\int_0^1 x^{\alpha-1}(1-x)^{2\beta-1}\left(1-\frac{x}{3}\right)^{2\alpha-1}\left(1-\frac{x}{4}\right)^{\beta-1} \,_a w_{p,b,c,\xi}^{\upsilon,\,\delta}\left(yx\left(1-\frac{x}{3}\right)^2\right)dx$$

$$= \left(\frac{2}{3}\right)^{2(\alpha+p+1)} \left(\frac{y}{2}\right)^{p+1} \Gamma(\beta) \,_2\psi_3 \left[ \begin{array}{c} (2\alpha+2p+2,\,4),\,(1,1); \\ \\ (\delta,\,\upsilon),\,\left(\frac{p}{\xi}+\frac{b}{2}+1,\,a\right),\,(2\alpha+\beta+2p+2,4); \end{array} \quad \frac{-4cy^2}{81} \right] \tag{2.5}$$

**Proof.** By making use of (1.1) in the ingrand of (2.4) and then interchanging the order of integral sign and summation, which is verified by uniform convergence of the involved series under the given conditions, we get

$$\int_0^1 x^{\alpha-1}(1-x)^{2\beta-1}\left(1-\frac{x}{3}\right)^{2\alpha-1}\left(1-\frac{x}{4}\right)^{\beta-1} \,_a w_{p,b,c,\xi}^{\upsilon,\,\delta}\left(y\left(1-\frac{x}{4}\right)(1-x)^2\right)dx$$

$$= \sum_{k=0}^{\infty} \frac{(-c)^k}{\Gamma(\upsilon k+\delta)\Gamma\left(ak+\frac{p}{\xi}+\frac{b+2}{2}\right)}\left(\frac{y}{2}\right)^{2k+p+1} \int_0^1 x^{\alpha-1}(1-x)^{2(\beta+p+2k)-1}\left(1-\frac{x}{3}\right)^{2\alpha-1}\left(1-\frac{x}{4}\right)^{\beta+p+1+2k-1} dx$$

$$\tag{2.6}$$

we can apply the integral formula (1.9) to the integral in (2.6) and obtain the following expression:

$$= \sum_{k=0}^{\infty} \frac{(-c)^k}{\Gamma(\upsilon k+\delta)\Gamma\left(ak+\frac{p}{\xi}+\frac{b+2}{2}\right)}\left(\frac{y}{2}\right)^{2k+p+1}\left(\frac{2}{3}\right)^{2\alpha}\frac{\Gamma(2\alpha)\Gamma(\beta+p+1+2k)}{\Gamma(2\alpha+\beta+p+1+2k)}$$

$$= \left(\frac{2}{3}\right)^{2\alpha}\left(\frac{y}{2}\right)^{p+1}\Gamma(2\alpha)\sum_{k=0}^{\infty}\frac{(-c)^k\Gamma(\beta+p+1+2k)}{\Gamma(\upsilon k+\delta)\Gamma\left(ak+\frac{p}{\xi}+\frac{b+2}{2}\right)\Gamma(2\alpha+\beta+p+1+2k)}\left(\frac{y}{2}\right)^{2k}$$

In accordance with the definition of (1.3), we obatain the result (2.5). This completes the proof of the theorem.

By similar manner as in proof of Theorem 3, we can prove the integral formula (2.6).

## 3. Special Cases

In this section, we derive some new integral formulae by using some known generalized struve function, which are given in collories 3.1 to 3.4.

 If we employ the same method as in getting Theorem 1 to 4, we obtain the following four corollaries with the help of (1.2) which is well known generalized Struve function due to Orhan and Yagmur [15, 16 ].

**Corollary 3.1.** *Let the condition of Theorem 1 be satisfied and for* $\upsilon = a = 1,\ \delta = 3/2$ *and* $\xi = 1,$ *theorem 1 reduces in following form*

$$\int_0^\infty x^{\mu-1}\left(x + a + \sqrt{x^2 + 2ax}\right)^{-\lambda} H_{p,b,c}\left(\frac{y}{x + a + \sqrt{x^2 + 2ax}}\right)dx$$

$$= 2^{-\mu-p}\, a^{\mu-\lambda-p-1}\, y^{p+1}\Gamma(2\mu)\ _3\psi_4\left[\begin{array}{c} (\lambda+p+2,\,2),\ (\lambda-\mu+p+1,\,2),\,(1,1); \\[2mm] \left(p+\dfrac{b+2}{2},\,1\right),\ (\lambda+p+1,\ 2),\,(\lambda+\mu+p+2,\ 2),\,(3/2,1); \end{array}\ \frac{-cy^2}{4a^2}\right] \quad (3.1)$$

**Corollary 3.2.** *Let the condition of Theorem 2 be satisfied and for* $\upsilon = a = 1,\ \delta = 3/2$ *and* $\xi = 1,$ *theorem 2 reduces in following form*

$$\int_0^\infty x^{\mu-1}\left(x + a + \sqrt{x^2 + 2ax}\right)^{-\lambda} H_{p,b,c}\left(\frac{xy}{x + a + \sqrt{x^2 + 2ax}}\right)dx$$

$$= 2^{-\mu-2p}\, a^{\mu-\lambda-1}\, y^{p+1}\Gamma(\lambda-\mu+1)\ _3\psi_4\left[\begin{array}{c} (\lambda+p+2,\,2),\ (2\mu+2p,\,4),\,(1,1); \\[2mm] \left(p+\dfrac{b+2}{2},\,1\right),\ (\lambda+p+1,\ 2),\,(\lambda+\mu+2p+2,\,4),\,(3/2,1); \end{array}\ \frac{-cy^2}{4}\right] \quad (3.2)$$

**Corollary 3.3.** *Let the condition of Theorem 3 be satisfied and for* $\upsilon = a = 1,\ \delta = 3/2$ *and* $\xi = 1,$ *theorem 3 reduces in following form*

$$\int_0^1 x^{\alpha-1}(1-x)^{2\beta-1}\left(1-\frac{x}{3}\right)^{2\alpha-1}\left(1-\frac{x}{4}\right)^{\beta-1} H_{p,b,c}\left(y\left(1-\frac{x}{4}\right)(1-x)^2\right)dx$$

$$= \left(\frac{2}{3}\right)^{2\alpha}\left(\frac{y}{2}\right)^{p+1}\Gamma(2\alpha)\ _2\psi_3\left[\begin{array}{c} (\beta+p+1,\,2),\,(1,1); \\[2mm] \left(p+\dfrac{b+2}{2},\,1\right),\ (2\alpha+\beta+p+1,2),\,(3/2,\,1),; \end{array}\ \frac{-cy^2}{4}\right] \quad (3.3)$$

**Corollary 3.4.** *Let the condition of Theorem 4 be satisfied and for* $\upsilon = a = 1,\ \delta = 3/2$ *and* $\xi = 1,$ *theorem 4 reduces in following form*

$$\int_0^1 x^{\alpha-1}(1-x)^{2\beta-1}\left(1-\frac{x}{3}\right)^{2\alpha-1}\left(1-\frac{x}{4}\right)^{\beta-1} H_{p,b,c}\left(yx\left(1-\frac{x}{3}\right)^2\right)dx$$

$$= \left(\frac{2}{3}\right)^{2(\alpha+p+1)}\left(\frac{y}{2}\right)^{p+1}\Gamma(\beta)\ _2\psi_3\left[\begin{array}{c} (2\alpha+2p+2,\,4),\,(1,1); \\[2mm] \left(p+\dfrac{b+2}{2},\,1\right),\ (2\alpha+\beta+2p+2,\,4),\,(3/2,1)\,; \end{array}\ \frac{-4cy^2}{81}\right] \quad (3.4)$$

**Conclusion:** Certain unified integral representation of generalized Struve function and its special cases are derived in this study. In this sequel, one can obtain integral representation of more generalized special function, which has much application in physics and engineering Science.

(1)Department of Mathematics, Wollo University, Dessie
P.O. Box: 1145, South Wollo, Amhara Region, (ETHIOPIA)
Email- dlsuthar@gmail.com

(2)Department of HEAS (Mathematics)
Rajasthan Technical University, Kota, (INDIA)
Email- sunil_a_purohit@yahoo.com

(3)Department of Mathematics
College of Arts and Science, Prince Sattam Bin Abdulaziz Unversity,
Wadi Al Dawaser, Riyadh Region, 11991, (SAUDI ARABIA)
Email- ksnisar1@gmail.com